\documentstyle[12pt]{article}

\title {The symplectic structure of the spin Calogero model.}
 \author{O. Babelon, M. Talon. 
 \thanks{L.P.T.H.E. Universit\'e Paris VI (CNRS URA 280), 
 Bo\^{\i}te 126, Tour 16, $1^{er}$ \'etage, 
 4 place Jussieu, F-75252 PARIS CEDEX 05} }
 \date{June 1997}

\begin{document}

\begin{titlepage}
\renewcommand{\thepage}{}
\maketitle
\vskip 2cm

\begin{abstract}
We compute the symplectic structure of the spin Calogero model in
terms of algebro-geometric data on the associated spectral curve.
\end{abstract}
\vfill
PAR LPTHE 97-30 \hfill
Work sponsored by CNRS: URA 280

\end{titlepage}
\renewcommand{\thepage}{\arabic{page}}

\section{Introduction}

In a recent paper~\cite{KrPh} I.M. Krichever and D.H. Phong made an
important progress in the understanding of the symplectic structure of
integrable models. Namely they were able to relate the canonical
symplectic structure of the model with a naturally defined
symplectic form on a suitable fibered space constructed with
algebro-geometric data. In their paper this construction is
illustrated by a number of examples. The next natural example to be
treated is the spin Calogero model, and we show in this paper that the
same general construction applies as well. This model has been solved
in~\cite{KBBT}, where it was found that as compared to the scalar case
a number of interesting new features enrich the algebro-geometric
analysis.

\section{The model.}

\subsection{Definition.}
The spin Calogero model consists of $N$ particles on a line
with internal degrees of freedom. The particles are described by their
positions $x_i$ and momenta $p_i$, together with spin variables
$f_{ij}$.
The Hamiltonian is:
\begin{equation}
H={1\over 2}\sum_i p_i^2-{1\over 2}\sum_{i\neq j}{f_{ij}f_{ji}\over
(x_i-x_j)^2}   \label{hami}
\end{equation}
The non vanishing Poisson brackets  for these  degrees of freedom are:
\begin{eqnarray}
\{p_i,x_j\}&=&\delta_{ij} \nonumber \\
\{f_{ij},f_{kl}\} &=& \delta_{jk} f_{il}- \delta_{il} f_{kj} \label{Kir1}
\end{eqnarray}
We  see  that the Poisson bracket~\ref{Kir1} is a Kirillov
bracket and therefore is degenerate. We shall choose an orbit such
that the matrix of spin variables $f$ is of rank $l$.

This system admits a Lax pair formulation, with Lax operator:
\begin{equation}
L_{ij}(t,z)=p_i \delta_{ij}+(1-\delta_{ij})f_{ij}\Phi(x_i-x_j,z)
\end{equation}
where for the elliptic model the special function $\Phi$ is defined
in terms of Weierstrass elliptic functions by:
$$\Phi(x,z)={\sigma(z-x)\over \sigma(z)\sigma(x)}e^{\zeta(z)x}$$
It is elliptic in $z$
and pseudo-periodic in $x$. It has been shown in~\cite{BAB} that the
model is integrable only when all $f_{ii}$ are equal, e.g. $f_{ii}=2$.

\subsection{The spectral curve.}

In the article~\cite{KBBT} it is noticed that this model can be
related to the matrix KP equation, and this leads to its solution by
algebro-geometric methods. In this paper we shall be concerned with
the Hamiltonian aspect of the solution. The main object of interest is
the spectral curve:
\begin{equation}
R(k,z)\equiv \det\,(2kI+L(t,z))=0 \label{spec}
\end{equation}
This defines a Riemann surface $\Gamma$ of genus $g=Nl-l(l+1)/2+1$
(see~\cite{KBBT}) which is time--independent due to the Lax equation.
The equation~\ref{spec} is polynomial in $k$ of degree $N$ with
coefficients elliptic functions of $z$, hence presents $\Gamma$ as an
$N$--sheeted covering of an elliptic curve with $2g-2$ branch
points. On such a curve one usually defines a line bundle by taking at
each point $P=(z,k)$ the eigenspace of $L(t,z)$ for the eigenvalue
$k$, i.e. the solutions of:
$$(2kI+L(t,z))C(t,P)=0$$

In particular at the branch points the eigenspaces for the two
colliding eigenvalues are generically of dimension 1, and $L(t,z)$ is
not diagonalizable for such $z$. 

Here special
attention is required for the points above $z=0$.

\subsection{The vicinity of $z=0$}

At $z=0$ the Lax matrix $L$ has an essential singularity since
$\zeta(z)= 1/z+O(z^3)$. Hence one can write:
$$L(z)={\rm Diag}(e^{{1\over z}(x_k-x_0)})\tilde{L}(z) {\rm
Diag}^{-1}(e^{{1\over z}(x_k-x_0)}) $$
where $x_0$ is some arbitrary origin, and $\tilde{L}$ is meromorphic
in a vicinity of $z=0$. The eigenvectors are of the form $C=(C_i)$
with: $$C_i(P)=e^{{1\over z}(x_i-x_0)}\tilde{C}_i(P)$$ where
$\tilde{C}_i(P)$ is a locally analytic eigenvector of $\tilde{L}$.
Since $\tilde{L}(z)=-{1\over z}(f-2I)+O(1)$ we see that $\tilde{C}_i(z=0)$
is an eigenvector of $f-2I$ and the corresponding eigenvalues are of
the form $k=(-1+\lambda/2)/z$, where $\lambda$ are the eigenvalues
of $f$, and $N-l$ of them vanish\footnote{
Remark that this defines $N$ generically different branches of
$\Gamma$ which intersect at the singular point $(z=0,k=\infty)$. After
one blow--up $l$ different points appear, and an ordinary multiple
point of order $N-l$ remains. An other blow--up at this point leaves
us with $N$ different points above $z=0$. When we speak of $\Gamma$ we
have in mind this desingularization, and we speak freely of the $N$
points $P_\alpha$ above $z=0$.}. 
The class of functions having
essential singularities at some points of the form 
$\psi=\exp(\alpha/z^m)\rho(z)$ (here $m=1$) with
$\rho$ locally meromorphic,  and meromorphic otherwise
are called Baker functions. The solution of many integrable models by
algebro--geometric methods essentially depends on the construction of
appropriate Baker functions on the spectral curve. A Baker function
has properties similar to ordinary meromorphic functions, e.g. has the
same number of zeroes and poles (consider the sum of residues of the
regular differential $d\psi/\psi$) but theorems such as Riemann--Roch
have to be modified, for example there exists non trivial Baker
functions with arbitrary prescribed $g$ poles on a surface of genus
$g$ while for a meromorphic function one needs to prescribe $g+1$
poles in the generic situation. To construct such a Baker function
$\psi$ let us remark that one can find a unique normalized abelian
differential of second kind $\omega_2$ and differential of third kind
$\omega_3( \rho_k)$ depending on the given poles, and the unknown zeroes
$\rho_k$ of $\psi$, such that $d\psi/\psi=\omega_2+
\omega_3( \rho_k)$. By definition this form integrates to 0 on
$A$--cycles. Imposing the same condition on $B$--cycles provides 
a system of $g$ equations for the $g$ unknowns $\rho_k$. For a construction using
theta functions see~\cite{KBBT}. Let us also remark that the quotient
of two Baker functions with the same type of singularities is a
meromorphic function to which one can apply the usual Riemann--Roch
theorem.

More generally one can define Baker line bundles, by imposing Baker
conditions on the components around singular points, and this is the
case for the eigenvector bundle of $L$.
This line bundle is of Chern class $-(g-1)$. 
In fact  in~\cite{KBBT} a nowhere vanishing
section $C(t,P)$ is constructed with $g-1$ poles. 

Since $L_{ij}$ is
not a symmetric matrix it is also convenient to construct an other Baker
line bundle by considering the adjoint eigenvector equation:
$$C^+(t,P)(2kI+L(t,z))=0$$
which is also of Chern class $-(g-1)$. 
Here the Baker behaviour occurs with opposite exponent:
$$C^+(P)=e^{{-1\over z}(x_i-x_0)}\tilde{C}^+(P)$$

Note that these Baker line
bundles are embedded in the ambient space of dimension $N$, which
provides a pairing:
$$<C^+(P),C(P')>=\sum_i C_i^+(P) C_i(P')$$
In particular  for two points $P$ and $P'$ above the
same $z$, corresponding to different eigenvalues, 
one has $<C^+(P),C(P')>=0$ due to the
eigenvector equations, hence at the branch points $<C^+(P),C(P)>$
vanishes for any regular sections $C^+$, $C$. Moreover the
singular factors $\exp {\pm 1\over z}(x_i -x_0)$ cancel above $z=0$ so
$<C^+(P),C(P)>$ extends to a meromorphic function on $\Gamma$.

Finally let $P_\alpha$ be the $N$ points above the same $z$, and  let
$C_\alpha = C(P_\alpha)$ be the corresponding eigenvectors. Except at the
branch points, they form a basis of ambient space, and any vector $V$
can be decomposed as:
\begin{equation}
V=\sum_{\alpha=1}^N {<C_\alpha^+,V>\over<C_\alpha^+,C_\alpha>}
C_\alpha \label{decomp}
\end{equation}

\subsection{Remarkable abelian forms.}

The Riemann surface $\Gamma$ presented by equation~\ref{spec}
possesses several remarkable abelian differentials. First the form
$dz$ is well defined on $\Gamma$ and has no pole. As such it has
$2g-2$ zeroes which are the branch points (where $k$ is the local
analytic parameter and $dz/dk$ vanishes). The form of main interest is
the form $kdz$ which has poles only above $z=0$. Finally let us take
two non vanishing sections $C^+$, $C$ of the above line bundles.
Then $<C^+(P),C(P)>$ is a meromorphic function without any
singularity above $z=0$ vanishing at the same points as $dz$ hence
$\Omega=dz/<C^+(P),C(P)>$ is an other analytic abelian form vanishing at the 
$2g-2$ poles $\gamma_k$ and $\gamma^+_k$ of $C$ and $C^+$. It will
play an important role later on.

\section{The canonical symplectic structure.}

Due to the degeneracy of the Poisson brackets~\ref{Kir1} one has to be
careful about the choice of the symplectic variety. In fact this
Poisson bracket is a Kirillov bracket for the coadjoint action of the
group $GL(N)$ acting on the dual of the Lie algebra $gl(N)$ identified
with itself by the invariant form $(A,B)\to {\rm Tr}(AB)$. The
orbits are generically characterized by the eigenvalues of the matrix
$f$ which are in the center of the Poisson bracket. Here we shall
consider matrices $f$ of rank $l$ with $l$ different non--vanishing
eigenvalues.  Such an orbit is of the form $\{g^{-1}\Lambda g|g\in
GL(N)\}$ with $\Lambda={\rm
Diag}(\lambda_1,\dots,\lambda_l,0,\dots,0)$.
The tangent space at the orbit of $f$ is the set of matrices $U=[f,X]$
for any $X\in gl(n)$. In a basis where $f$ is
diagonal this equation reads $U_{ij}=(\lambda_i -\lambda_j)X_{ij}$, hence $U_{ij}$
vanishes when $\lambda_i =\lambda_j$ but is otherwise arbitrary. So
the dimension of the orbit is $2Nl-l^2 -l$. At a point $f$ the
 symplectic form on two tangent vectors
$[f,X]$ and $[f,Y]$ is given by the well-defined formula:
$$\omega_K([f,X],[f,Y])={\rm Tr}(f[X,Y])$$
in accordance with Kirillov's prescription. Explicitly in a basis
where $f$ is diagonal this can be written:
\begin{equation}
\omega_K = \sum_{{i,j\atop\lambda_i\neq\lambda_j}} {df_{ij}\wedge df_{ji}
\over \lambda_i-\lambda_j} \label{Kir2}
\end{equation}
where by definition on a tangent vector $U$, $df_{ij}(U)=U_{ij}$.

Of course the symplectic form of our model is 
\begin{equation}
\omega=\sum dp_i\wedge
dx_i+ \omega_K  \label{canonique}
\end{equation}

The Hamiltonian~\ref{hami} is not invariant under the
above $GL(N)$ but it is is preserved by special subgroups. First we
have the discrete subgroup of permutation matrices, i.e. the Weyl
group, which simply operates by permutation of the $N$ indices $i$
\footnote{These actions extend obviously to actions on $L$ by
similarity transformations},
and more importantly we have the group of diagonal matrices, i.e. the
Cartan torus, which operates by:
$$f_{ij}\to d_i^{-1}f_{ij}d_j$$ This action preserves the Hamiltonian which
 only depends on $f_{ij}f_{ji}$ and all higher Hamiltonians ${\rm Tr}(L^n)$.
The action of this toral subgroup induces a fibering of the
orbits into fibers of dimension $N-1$ since multiple of the identity
leave $f$ invariant. The moment associated to this action is the
collection of diagonal elements $f_{ii}$, that is $N-1$ non trivial
moments since on the orbit the eigenvalues of $f$ are fixed, hence so
is ${\rm Tr}(f)$. We consider the reduced dynamical system obtained by
first fixing the moments to a common value $f_{ii}=2$, and then
quotienting by the stabilizer of this moment which is the whole
diagonal group. This is the integrable system as considered above.

We can now count the number of degrees of freedom. We have $2N$
degrees for the $x_i, p_i$, plus $2Nl-l^2 -l$ for the orbit, minus
$2(N-1)$ due to the Hamiltonian reduction which ends up to a phase
space of dimension $2(Nl-l(l+1)/2+1)=2g$.

It is a remarkable fact that the spectral equation~\ref{spec} is
dependent on $g$ non trivial integrals of motion. In fact it has been
shown in~\cite{KBBT} that considering only the order of the singularity at $z=0$
and knowing that all $f_{ii}=2$ and $f$ is of rank $l$, the spectral
equation depends on $g+l-1$ parameters. Here the eigenvalues
of $f$ are in the center of the Poisson algebra, and must not be
counted as dynamical variables. There are
 $l-1$ independent non vanishing eigenvalues of $f$
 since Tr$(f)$ has previously
been fixed to $2N$. Thus we end up with exactly $g$ action variables.

Our task is to find $g$ other dynamical variables which will be used
to construct the angle variables. It is known that such algebraically
integrable systems linearize on the Jacobian of the spectral
curve, hence it is natural to use $g$ points $\gamma_k$ on $\Gamma$ as
complementary variables. We have  found above $g-1$  poles of a
non--vanishing section of the eigenvector bundle.  If one defines
$C(P)$ such that $C_1(P)=1$ as in ~\cite{KBBT}, the poles of $C$ are
given by the vanishing on $\Gamma$ of the first minor of the matrix
$L+2kI$. This is the algebraic equation which relates the $\gamma_k$
to the dynamical variables of the system. Note that in this minor the
variables $x_1$ and $p_1$ have disappeared. This corresponds to
reduction by translational symmetry, which leaves a phase space of
dimension $2(g-1)$. 

In order to get the full phase space it is
convenient to take the product
of such a section with an appropriate Baker function yielding a
section with $g$ poles and a fixed zero. It suffices to choose it with
given zeroes canceling the poles of $C$, and one more zero 
at a fixed point $P_0$ above
$z=0$. The singular behaviour above $z=0$ is taken to be $\exp
(x_0/z)$ so that we end up with $C_i(P)$ proportional to $\exp (x_i/z)$
having $g$ dynamical poles, and similarly for $C^+$. Note that the
poles $\gamma_k^+$ of $C^+$ are determined when the poles $\gamma_k$
of $C$ are given, since the abelian form $\Omega=dz/<C^+(P),C(P)>$
is meromorphic with only singularity a double pole at $P_0$, hence
has no residue at $P_0$, and one can choose the global normalization
of $C^+$ such that:
$$ \Omega = ({1\over z^2}+O(1))dz \quad P\to P_0$$
But such a form is uniquely determined when one  fixes $g$ zeroes
$\gamma_k$, and it has $g$ other determined zeroes $\gamma_k^+$ (otherwise the 
quotient of two such forms would be a meromorphic function with
$g$ poles and $g$ other zeroes, which is generically forbidden).
We summarize the definition of $C$ and $C^+$ by stating their
behaviour as $z\to 0$:
\begin{equation}
C_i(P)=e^{x_i\over z} c_i(P), \quad
C_i^+(P)=e^{-x_i\over z} c_i^+(P), \quad z\to 0
\label{singul}
\end{equation}
with $c_i$ and $c_i^+$ regular (and vanish when $P\to P_0$).
This together with the above condition on $\Omega$ clearly
determines $C$ and $C^+$ up to a constant factor $\lambda$
on $C$ and $1/\lambda$ on $C^+$ when the $\gamma_k$ are given,
i.e. when the dynamical variables are given. Moreover, the
$c_i(P)$ for the $N$ points $P$ above $z=0$ are 
$N$ eigenvectors of $f$ (for $P=P_0$ of course we take
$c_i(P)/z$ as the corresponding eigenvector), and similarly
for $c_i^+$.

\section{The action--angle variables}

\subsection{Some fiber bundles.}

Following the ideas of~\cite{KrPh} we first introduce some natural bundles on the
moduli space of our curves, or more specifically the $g$--dimensional space
of action variables. Let us stress that we keep the eigenvalues of $f$
fixed throughout our discussion. First we have the bundle $\cal G$ whose fiber
above a  particular spectral function $R(k,z)$ is the curve $\Gamma$ of
equation $R=0$. Note that the differential $dz$ is naturally defined on the 
fibers independent of $R$. We denote $\delta$ the differential on $\cal G$,
and we see that it can be splitted into a vertical part along the fiber
$\delta_V$ and an horizontal part which corresponds to differentiation
with $z$ fixed $\delta_H$. Hence if $a_i,\; i=1,\dots, g$ are independent
action variables we can write:
$$\delta=\sum_i {\partial\over\partial a_i}da_i + {\partial\over\partial z}dz $$
In particular the two--form $\delta (kdz)=\delta_H k \wedge dz$ is well defined
on $\cal G$. Note that $k$ has simple poles on $\Gamma$ whose residues are
given by the eigenvalues of $f$ as seen above, hence the horizontal
differentiation of these residues vanish, so the above form is regular
on $\cal G$. This is one of the key observations in~\cite{KrPh}.

Similarly, we define the fiber bundle $\cal J$ whose fiber above $R$
is the Jacobian of the curve of equation $R=0$. A point of the Jacobian 
can be seen generically as
an unordered set of $g$ points $\gamma_i,\; i=1,\dots,g$ of the curve
$\Gamma$, and we also denote $\delta$ the differential on $\cal J$
which can be splitted as above, simply replacing the vertical part by:
$$\delta_V=\sum_i {\partial\over\partial z_i}dz_i$$
where $z_i$ is locally the $z$--coordinate of $\gamma_i$.
It follows that we can define~\cite{KrPh} on $\cal J$ a regular two--form
$\omega=\delta (\sum_i k_i dz_i )$ which defines a symplectic structure.
The notation is not accidental since we shall see that $\omega$
precisely reduces to the canonical symplectic structure when we identify
$\cal J$ to the phase space of our dynamical system.

It is important to define, as in~\cite{KrPh}, forms on $\cal J$ with
values abelian forms on the curve $\Gamma$ associated to
the corresponding base point of moduli space. In fact these
forms are defined on the bundle whose fiber above a base point
$R$ is the Cartesian product of the curve $\Gamma$ of equation $R=0$
and its Jacobian Jac$(\Gamma)$. A generic point of such a bundle
is described by $g$ action variables $a_i$, a point $P$ on
$\Gamma$, and a set of $g$ points $\gamma_k$ of $\Gamma$. 
The differential $\delta$ acts on the $a_i$ and the $\gamma_k$,
while the differential $d$ acts vertically along $\Gamma$ on $P$.

We have found particularly
convenient to introduce the meromorphic two--form on $\cal J$:
\begin{equation}
\Phi = < \delta(C^+(P) \Omega(P))\wedge \delta L\, C(P) > 
\label{astuce}
\end{equation}
which due to the $P$ dependence and the presence of $\Omega(P)=
dz/<C^+(P)C(P)>$ has values one--form on $\Gamma$. Here $L$
is the Lax matrix depending implicitly on the $a_i$ and $\gamma_k$,
and the phase space differentiations appear in
$\delta (C^+\Omega)$ and $\delta L$. Finally the brackets $<>$ contract
the $i$--indices of $C$, $L$ and $C^+$. The main trick~\footnote{
A similar trick allows a considerable simplification of the proof of the main
theorem in~\cite{KrPh}.} is to consider the
sum of the residues of $\Phi$ on $\Gamma$ which must vanish. But
$C^+(P)\Omega(P)$ is regular at $\gamma_k^+$ for any value of the moduli
and the $\gamma_k$ hence so is its $\delta$--differential. Hence the
only possible singularities of $\Phi$ are located at the $g$ points $\gamma_k$ of
$\Gamma$ where $C(P)$ has a pole,
and the $N$ points above $z=0$ where $\Omega$ and $\delta L$ are singular.

Writing these residues, we get a relation between
a two--form on $\cal J$ which happens to be the previously defined
algebro--geometric 
symplectic form, and a two--form located above $z=0$ which
boils down to the canonical symplectic form of our model.

\subsection{Algebro--geometric description of the canonical
symplectic structure.}

We can now state the result of this paper:
\proclaim Proposition.
1) The canonical symplectic form~\ref{canonique} of the spin Calogero model can
be written in terms of algebro--geometric variables as:
\begin{equation}
\omega = \sum_{i=1}^g \delta k_i \wedge \delta z_i \label{skly}
\end{equation}
where $z_i=z(\gamma_i),\; k_i=k(\gamma_i)$. \hfill\break\noindent
2) If we take as angle variables the Abel transform of the divisor
$\sum \gamma_i$ namely the angles defined modulo periods of $\Gamma$:
\begin{equation}
\theta_k=\sum_i \int^{\gamma_i} \omega_k  \label{abel}
\end{equation}
where the $\omega_k$ are a basis of regular abelian differentials
dual to a basis $\{A_k\}$ of $A$--cycles  of $\Gamma$, then the canonically
conjugated variables are given by:
\begin{equation}
a_l=\oint_{A_l} k dz \label{canon}
\end{equation}

Several remarks are in order at this point. First let us note that
the first statement is exactly the same type of result advocated
by Sklyanin in his famous solutions of various integrable models
by his method of separation of variables~\cite{Skl}. In fact when
one takes into account the above description of the poles of $C$
one sees that it is really the same statement. However the result
is obtained here without any appeal to the $R$--matrix method.
Moreover the second statement shows independently that the motion
linearizes on the Jacobian of the spectral curve.

Then let us note that formulae of type~\ref{canon} have already
appeared at various places. For example it is clear that they play
a role in the classical and semi--classical analysis of the Neumann
model~\cite{BT}, and have been introduced more generally in ~\cite{FlMl}. In 
more recent analysis of quantum integrable models~\cite{Smi} deformations
of these formulae play a central role. One may hope that direct
deformation of this description leads to the quantification of integrable
systems in terms of algebro--geometric concepts.

Finally, note that the proof of $1) \Rightarrow 2)$ is straightforward
using a clever argument of~\cite{KrPh}. In fact let us define $a_i$
according to equation~\ref{canon}, they are obviously moduli
coordinates on the basis of our fiber space. Now as noted above
$\delta(kdz)$  defined on $\cal G$ is an analytic form on the fibers,
hence can be expanded on the basis $\omega_k$ of analytic one--forms.
To find the coefficients, let us compute:
$$\oint_{A_j} {\partial (kdz)\over \partial a_i}={
\partial a_j \over \partial a_i}=\delta_{ij}$$
so we see that ${\partial (kdz) \over \partial a_k}=\omega_k$ and we have:
$$\delta (kdz)=\sum_i \delta a_i \wedge \omega_i$$
Now, taking into account the value~\ref{abel} of the angular variables
we write the symplectic form as:
$$\omega=\sum_j \delta(k_j \delta z_j)=\sum_{ij} \delta a_i \wedge\omega_i(\gamma_j)
=\sum_i \delta a_i \wedge \delta \theta_i$$
which shows that the $a_i$ are indeed canonically conjugated to the $\theta_i$.
One should note the simplicity of this derivation which
essentially uses only one ingredient: the fact that the residues of the poles
of $kdz$ are killed by $\delta$ so that $\delta(kdz)$ is regular.
This can be contrasted to the involved computations inherent to previous
approaches on this subject.

\subsection{Proof of the main result.}

Let us compute the sum of the residues of the one--form $\Phi$ on
$\Gamma$ (which are themselves two--forms on $\cal J$). As noted above
we have to look at the points $\gamma_k$'s and at the points
above $z=0$. In the vicinity
of a pole $\gamma_k$ we can write $C^+ \Omega=(z-z_k)\Psi$ with $\Psi$
regular for $z=z_k$, hence:
$$\delta(C^+ \Omega)=-C^+\Omega {\delta z_k \over z-z_k}+{\rm
regular}$$
(here $\delta z_k=\delta_V z_k$), so the residue at $\gamma_k$ is:
$${\rm Res}_{\gamma_k} \Phi=-\delta_V z_k \wedge {<C^+ \delta L\, C>
\over <C^+C>}(\gamma_k)$$
Note that $\delta_H (C^+\Omega)$ does not contribute to this residue
because by definition $\delta_H z_k=0$.
But $C^+ L=-2kC^+$ and $LC=-2kC$ hence $\delta L C+L\delta C=
-2k \delta C -2\delta k C$ which  upon bracketing with $C^+$ yields
$<C^+ \delta L\, C>=-2\delta k <C^+C>$ so the final expression for the
sum of the residues at the $\gamma_j$ is:
$$\sum {\rm Res}_{\gamma_j}= - 2 \sum_j \delta k_j \wedge \delta z_j
= - 2 \omega$$
where $\omega$ is the previously defined symplectic form on $\cal J$.
Writing that the sum of the residues vanishes we have:
\begin{equation}
2 \omega= \sum_\alpha {\rm Res}_{P_\alpha} \Phi \label{Main}
\end{equation}
where the $P_\alpha$'s are the $N$ points above $z=0$.

Our next task is to compute this sum and show that it indeed reduces
to the canonical symplectic form on our dynamical system.
Let us remark that in the vicinity of $z=0$ we have a common local
parameter $z$ for all the curves in $\cal G$, and $\Omega$ is defined
by the same normalization at $P_0$ so $\delta \Omega$ is regular
at $P_0$.
Due to the normalization~\ref{singul} we can write:
$$\delta C^+={-x_i\over z}C^+ + e^{-x_i\over z}\delta c^+$$
We treat separately the contribution of the first term.
The considered expression in the vicinity of $P_\alpha$ takes the
form:
$$-\sum_{i,j} {dz\over z} {<C^+_i \delta x_i\wedge \delta L_{ij} C_j>
\over <C^+C>}(P)$$
Let us introduce the matrix $M_{ij}=\delta x_i\wedge \delta L_{ij}$
which only depends on $z(P)$. Summing on the $N$ sheets we get, in
view of equation~\ref{decomp}\footnote{
The point $P_0$ presents no special problem since the factor $z^2$
present in the numerator cancels a similar factor in the denominator.}, 
the simple expression $-dz/z {\rm
Tr}(M)$, so the residue at $z=0$ is just $\sum_i \delta p_i \wedge \delta x_i$.
For the same reason the contribution proportional to $\delta \Omega$
can only pick diagonal terms in $L$ and there is no residue.

So we are left with the remaining contribution of $e^{-x_i\over z}\delta c^+$.
Here we must examine the $1/z$ terms present in $\delta L$. First note
that all factors $e^{\pm x_i\over z}$ cancel between $e^{-x_i\over z}\delta c^+$,
$\delta L_{ij}$ and $C_j$. Then taking into account the identity:
$${\partial\over\partial x}\phi(x,z)=\phi(x,z)[\zeta(x+z)-\zeta(x)]$$
and the expansion $\phi(x,z)=(-1/z+\zeta(x)+O(z))\exp(\zeta(z)x)$ one
sees that all terms in $\delta p_i$ and $\delta x_i$ are regular at
$z=0$, so only terms in $\delta f_{ij}$ contribute to the residue.
Moreover $(1-\delta_{ij})\delta f_{ij}=\delta f_{ij}$ because
$\delta f_{ii}=0$ so one can simply replace $\delta L_{ij}$ by
$(-1/z)\delta f_{ij}$. The residue at $z=0$ is now obvious and
the expression to be computed is therefore:
$$-\sum_\alpha {<\delta c^+_\alpha \wedge \delta f  c_\alpha>
\over <c^+_\alpha c_\alpha>}$$
where $c_\alpha$ is the value of $c(P)$ at $P_\alpha$, i.e. they are
the eigenvectors of $f$, and similarly for $c^+$.

To complete the calculation note that varying the equation $fc_\alpha=
\lambda_\alpha c_\alpha$ keeping $\lambda$ constant 
(recall that eigenvalues of $f$ are fixed) one gets by contracting
with $c^+_\beta$ the relation $<c^+_\beta \delta f c_\alpha>=
(\lambda_\alpha - \lambda_\beta)<c^+_\beta \delta c_\alpha>=
- (\lambda_\alpha - \lambda_\beta)<\delta c^+_\beta  c_\alpha>$.
Expanding the above residue on the basis $c_\beta$ in view
of~\ref{decomp} one gets:
$$-\sum_{\alpha\beta} {<\delta c^+_\alpha  c_\beta>\over <c^+_\beta  c_\beta>}
\wedge {<c^+_\beta \delta f c_\beta>\over <c^+_\alpha  c_\alpha>}$$ in which
the sum can be restricted to $\lambda_\alpha \neq \lambda_\beta$ by
the previous relation. But using~\ref{decomp} this is
nothing more than:
$$+\sum_{\alpha,\beta\atop \lambda_\alpha \neq \lambda_\beta}
{\delta f_{\alpha\beta}\wedge \delta f_{\beta\alpha} \over
\lambda_\alpha - \lambda_\beta}$$
in which we recognize the expression~\ref{Kir2} of the Kirillov form
$\omega_K$. Finally we have proven that:
$$2\omega=\sum_i \delta p_i \wedge \delta x_i+\omega_K$$

\section{Conclusion.}

The example treated in this paper shows once more the power
and the generality of the method introduced by
Krichever and Phong. Even in an intricate dynamical situation, the
canonical symplectic form can be expressed in terms of
algebro--geometric data with the same formula as in standard cases.
Due to its simplicity and generality, this formula may very well
extend to the quantum domain, as simple examples already indicate.

{\bf Acknowledgements.} We thank I. Krichever and F. Smirnov for
useful discussions.


\begin{thebibliography}{xxxxxxxx}


\bibitem{KrPh} I.M. Krichever and D.H. Phong, {\it
On the integrable geometry  of soliton equations and
the $N=2$ supersymmetric gauge theories.} HEP-TH 9604199.

\bibitem{KBBT} I.M. Krichever, O. Babelon, E. Billey and M. Talon,
{\it
Spin generalization of the Calogero-Moser system and the matrix KP
equation.}
Amer. Math. Soc. Transl. (2) {\bf 170} 1995  83--119, HEP-TH 9411160.

\bibitem{BAB} E. Billey, J. Avan and O. Babelon, {\it Exact Yangian
symmetry in the classical Euler-Calogero-Moser model.}
Phys. Lett. A {\bf 188} (1994) 263--271 HEP-TH 9401117.

\bibitem{Skl} E.K. Sklyanin, {\it Separation of variables. New
trends.}
Progr. Theor. Phys. Suppl. {\bf 118} (1995) 35--60.

\bibitem{BT} O. Babelon and M. Talon, {\it Separation of variables
for the classical and quantum Neumann model.}
Nucl. Phys. {\bf B379} (1992) 321--339 HEP-TH 9201035.

\bibitem{FlMl} H. Flashka and D.W. McLaughlin, Progr. Theor. Phys.
{\bf 55} (1976) 438--456.

\bibitem{Smi} F.A. Smirnov, {\it On the deformation of abelian
integrals.} Lett. Math. Phys. {\bf 36} (1996) 267--275 Q-ALG 9501001


\end{thebibliography}
\end{document}